\documentclass[12pt]{article}

\usepackage{amsfonts,amsmath,graphicx}
\def\id{\mathop{\rm id}}

\newtheorem{lemma}{Lemma}[section]
\newtheorem{theorem}{Theorem}[section]

\newtheorem{corollary}{Corollary}[section]

\title {The Bennequin Number, Kauffman Polynomial, and Ruling Invariants of A 
Legendrian Link:\\ the Fuchs Conjecture and Beyond}

\author{Dan RUTHERFORD}

\begin{document}

\maketitle

\begin{abstract} We show that the ungraded ruling invariants of a legendrian link can be realized as certain coefficients of the Kauffman polynomial which are non-vanishing if and only if the upper bound for the Bennequin number given by the Kauffman polynomial is sharp.  This resolves positively a conjecture of Fuchs.  Using similar methods a result involving the upper bound given by the HOMFLY polynomial and 2-graded rulings is proved.
\end{abstract}

\section{Introduction}

 Historically, the first examples of invariants distinguishing between 
Legendrian links of the same topological type were given by the {\it Bennequin 
number}, $\beta (K)$, and the {\it rotation number} $r(K)$.  The Bennequin 
number can be negative of arbitrary large magnitude within a topological link 
class, but perhaps the first substantial result of the theory was Bennequin's 
upper bound,$$\beta(K) + |r(K)| \leq 2g(K) -1$$where $g(K)$ is the genus of 
$K$ (see \cite{B}).

This result instigated much future research which I divide roughly into 2 
directions.

1.  Find new upper bounds for $\beta(K) $ in terms of topological link 
invariants, and clarify the relationship between different bounds for 
$\beta(K)$.

2.  Determine the maximal value for $\beta(K)$ within a given knot type, or 
more generally characterize classes of knots for which a given bound is sharp.

Many results exist in the direction of 1., see \cite{Ng1} for a list of 
currently known bounds, most of which are known to be independent of one 
another.  In the direction of 2., at present, the maximal value of $\beta(K)$ 
within a given topological class has been tabulated for knots of 9 crossings 
and less.  Several partial results also exist clarifying when specific bounds 
for $\beta(K)$ are sharp (see \cite{EF} cite{Ng1} \cite{Ng2} \cite{T}).

The main results of this paper give simple necessary and sufficient 
conditions for the upper bounds on $\beta(K)$ given by Kauffman polynomial, 
$F_K(z,a)$, $$\beta(K) <  -\deg_aF_K$$and by the HOMFLY polynomial $P_K(z,a)$,
$$\beta(K) + |r(K)|< -\deg_aP_K(z,a)$$to be sharp (as estimates for $\beta 
(K)$). 

In the case of the Kauffman polynomial bound, the equivalent condition was 
conjectured by Fuchs in \cite{F}, and is precisely the existence of an 
(ungraded) ruling of a front diagram for $K$.  Interestingly enough, this 
condition is itself known to be equivalent to the existence of an augmentation 
on the Legendrian contact DGA \cite{F} , \cite{FI}, \cite{S}  as defined by 
Chekanov and Eliashberg  \cite{Ch},\cite{El2}.  The number of ($p$-graded) 
rulings with the $\#\{$switches$\} - \#\{$left cusps$\} =n$ fixed (see Section 
2) has been considered as a combinatorial invariant \cite{ChP} , with the 
sequence over all $n$ referred to as the complete ruling invariant in 
\cite{NgS}.  In the case of $p$-graded rulings these invariants can 
distinguish knots with identical classical invariants.  On the contrary, our 
result shows that in the ungraded case the complete ruling invariant is given 
by certain coefficients of $F_K$ depending only on $\beta(K)$, hence depends 
only on topological knot type and classical invariants (as conjectured in 
\cite{NgS}).  For knots, the 2-graded ruling invariant is realized in 
coefficients of the HOMFLY polynomial. 

Together these results clarify the relationship between these two upper bounds 
in an interesting way.  It is known \cite{Fer} that no inequality exists in 
general between $\deg_aP_K$ and $\deg_aF_K$.  However, in the case when the 
HOMFLY estimate is sharp, the Kauffman estimate must be sharp as well. 
Finally, together with a proposition of Ng \cite{Ng1}, our result implies that 
the Kauffman bound is sharp for alternating links.  As Ng notes, this 
positively answers a question of Ferrand \cite{Fer}, who asked whether for 
alternating links the estimate coming from the Kauffman polynomial should be 
better than that given by the HOMFLY polynomial.

\subsection{Acknowledgements} I would like to express my deep thanks to Dmitry 
Fuchs for introducing me to both the subject and his conjecture, and for 
creating the figures which appear in this document.  Also, thank you to Chris 
Berg.  In the summer of 2004 we showed that the Kauffman estimate was sharp 
for many families of Legendrian links which admit rulings but did not publish 
our results.   This research was supported in part by NSF VIGRE Grant No. DMS-0135345.

\section{Preliminaries and Definitions}

\subsection{Bennequin numbers and Kauffman polynomials}

\begin{figure}[hbtp]
\centering
\includegraphics[width=3.6in]{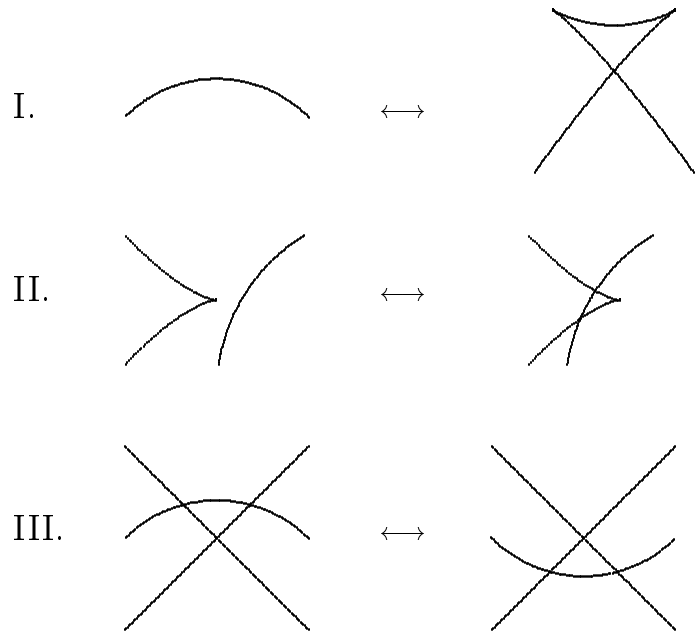}
\caption{Legendrian Reidemeister moves}
\end{figure}

We study Legendrian knots and links in ${\bf R}^3$ with respect to the 
standard contact structure given by the kernel of the $1$-form $y dx - dz$.  
We adopt a diagrammatic perspective where knots are presented by their 
projection into the $xz-$plane, hence forth referred to as the {\sl front 
diagram} or {\sl front}.  A front diagram is a smooth map of (a disjoint union 
of several copies of) $S^1$ into the plane with no vertical tangents and 
no singularities save cusps and non-tangential double points.  Two Legendrian 
links are {\sl Legendrian isotopic} if there fronts can be transformed into 
one another through a sequence of the following three {\sl Legendrian 
Reidemeister moves} (see Figure 1) and planar isotopies through front 
diagrams.  

Given a front diagram of a link $K$ by smoothing cusps and placing the strand 
with lesser slope on top at crossings we arrive at a diagram of an ordinary
topological knot class.  We denote this diagram as $Top(K)$.  All topological 
invariants of $Top(K)$ will form Legendrian invariants of $K$, since
Legendrian equivalence is strictly stronger.

Given an oriented front $K$ let $c(K)$, $cr(K)$ and $w(K)$ be the number of 
left cusps of K, the number of crossings of K, and the writhe of $Top(K)$ 
respectively.  The Bennequin number is defined by \begin{equation} \beta(K) = 
w(K) - c(K).\end{equation}

For knots $\beta(K)$ is independent of orientation.  The Bennequin number can 
easily be decreased within a topological link type by adding {\sl zig-zags} 
(see Figure 2).  

\begin{figure}[hbtp]
\centering
\includegraphics[width=3.4in]{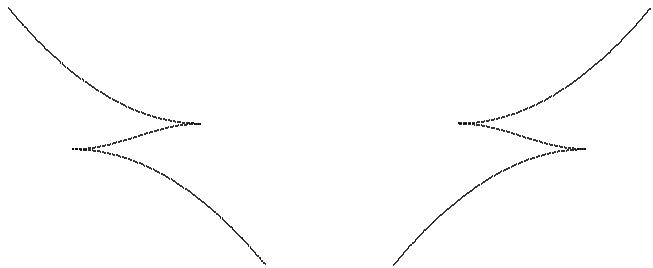}
\caption{Zig-zags}
\end{figure}

The following is one of many estimates showing that $\beta(K)$ is 
bounded above within a topological link type.  

\begin{lemma}[\cite{FT}, \cite{Fer}]

$\beta(K) <  -\deg_aF_K$ where $F$ denotes the two variable 
Kauffman Polynomial.

\end{lemma}

To define the Kauffman Polynomial (Dubrovnik version) \cite{K1} first an 
auxiliary polynomial, $D_K(z,a)$, is defined up to regular isotopy.  $D_K$ is 
characterized by skein relations 

\begin{figure}[hbtp]
\centering
\includegraphics[width=3in]{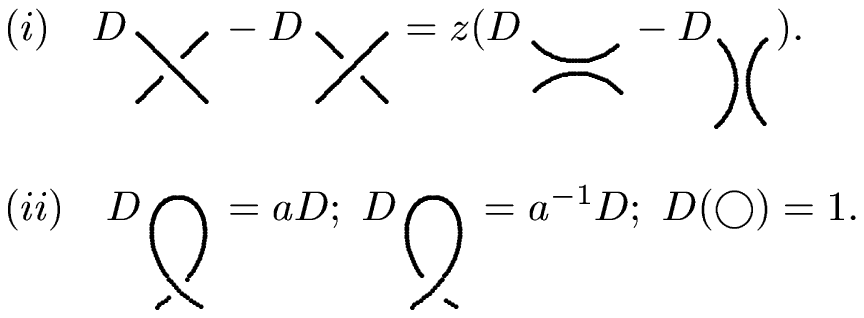}
\end{figure}

\noindent and invariance under type II and type III Reidermeister moves. The 
Kauffman polynomial $F_K$ is then defined as $F_K := a^{-w(K)} D_K$. $D_K$ 
depends on the diagram representing $K$, so if we are given a front $K$ we let 
$D_K := D_{Top(K)}$.

Lemma 2.1 has a simpler equivalent statement in terms of the $D$ polynomial.

\begin{lemma}

For any front $K$, $c(K) - 1 \geq \deg_aD_K$ with equality if and only if 
$\beta(K) =  -\deg_aF_K -1$

\end{lemma}

\paragraph{Proof.}  \[\begin{array} {l} w(K) - c(K) =\beta(K)   \leq -
\deg_aF_K -1 = w(K) -\deg_aD_K -1\\  \Leftrightarrow c(K) -1 \geq \deg_aF_K.
\end{array}\]

Note that Lemma 2.2 implies that the sharpness of Lemma 2.1 for a given 
Legendrian link class is independent of orientation.

\subsection{Rulings and $R_K$}

The following notion of a {\sl ruling} was introduced independently (but for 
entirely different purposes) by Fuchs \cite{F} and Chekanov and Pushkar 
\cite{ChP}. A similar notion (but without the normality condition (iii)) was 
considered as early as in 1987 by Eliashberg \cite{El2}.

By planar isotopy we may assume that all singularities of a front $K$ have 
different $x$-coordinates, and we will do so henceforth.  Given a subset $\rho 
=\{\lambda_1, \ldots \lambda_M\}$ of the set of crossings of $K$,  with the 
$x$-coordinate of $\lambda_i$ denoted $x_i$ so that $x_i < x_{i+1}$, let 
$S_\rho(K)$ denote the front obtained from $K$ by resolving all crossings in 
$\rho$ to parallel horizontal lines (see Figure 3). 

\begin{figure}[hbtp]
\centering
\includegraphics[width=5in]{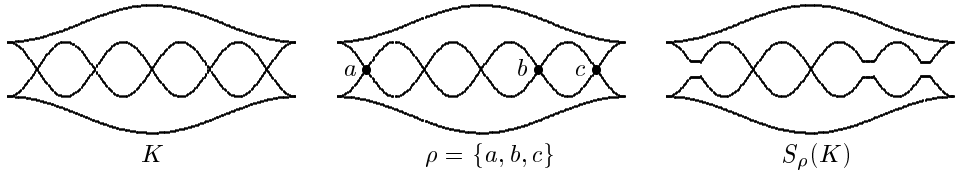}
\caption{A ruling}
\end{figure}

The set $\rho$ is called a ruling if

(i)  every component $T_j$ of $S_\rho(K)$ consists of two horizontal strands 
having one left cusp and no self crossings.  The upper is denoted $U_j$, and 
the lower $L_j$. 

(ii) for each $i$, the strands of $S_\rho(K)$ meeting where $\lambda_i$ was in 
$K$ belong to different components.  Call the upper of these strands $P_i$ and 
the lower $Q_i$.

(iii) one of the following normality conditions holds for each $i$:

For some $j_1, j_2$,

(a) $P_i = L_{j_1}$ and $Q_i = U_{j_2}$ 

(b) $P_i = U_{j_1}$ and $Q_i = U_{j_2}$, with the $z$-coordinate of $L_{j_1}$ 
less than the $z$-coordinate of $L_{j_2}$ at $x = x_i$

(c) $P_i = L_{j_1}$ and $Q_i = L_{j_2}$, with the $z$-coordinate of $U_{j_1}$ 
less than the $z$-coordinate of $U_{j_2}$ at $x = x_i$.(See Figure 4.)

\begin{figure}[hbtp]
\centering
\includegraphics[width=5in]{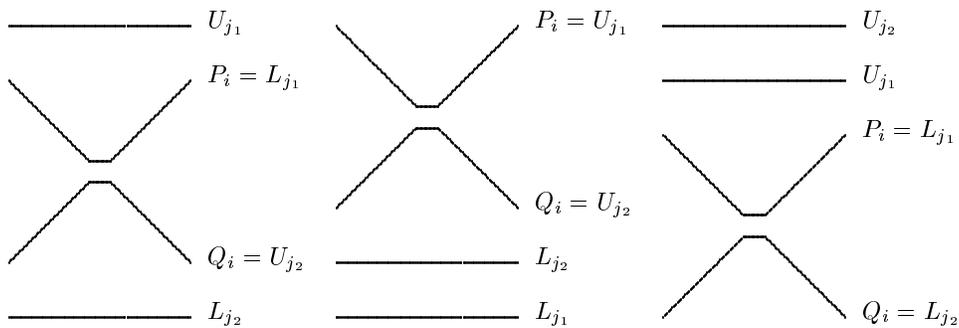}
\caption{Normality condition}
\end{figure}

For example, the ruling shown in Figure 3 meets the normality condition, but 
it will fail to do so, if we shift $b$ to the next crossing at the left.

The elements of a ruling $\rho$ are called {\sl switches} and we denote the 
number of switches in $\rho$ as $s(\rho) := \#\rho $.

For a front $K$, let $\Gamma (K)$ be the set of rulings of $K$. 

For integer $n$, let $f_n:= \# \{\rho \in \Gamma(K) | s(\rho) - c(K) +1 = 
n\}$.  To simplify future notation set $$j(\rho) := s(\rho) - c(K) +1$$ The 
sequence $\{f_n\}$ is a Legendrian Isotopy invariant \cite{ChP},(this can 
easily be seen by constructing bijections between rulings under Legendrian 
Reidemeister moves) and we condense it into the {\it Ruling Polynomial},

\begin{equation} 
R_K := \sum_{n \in {\bf Z}}f_nz^n = \sum_{\rho \in \Gamma(K)} z^{j(\rho)} 
\end{equation}

\subsection{Definition of the polynomial $B$.}

Let $D_K(z,a) = \sum_{n \in {\bf Z}} C_n(z)a^n$.  Define $$B_K(z) = C_{c(K) 
-1}(z)=\ \mbox{coefficient\ of}\ a^{-1}\ \mbox{in}\ a^{\beta(K)}F_K.$$

$B_K$ is a Legendrian isotopy invariant and, by Lemma 2, is non-zero iff the 
Kauffman estimate for the Bennequin number is sharp.

\section{Main result}

\begin{theorem}
For any Legendrian link $K$, $R_K = B_K$.
\end{theorem}
\begin{lemma}
$R_K$ and $B_K$ both satisfy the following skein relations:
\end{lemma}

\begin{figure}[hbtp]
\centering
\includegraphics[width=4in]{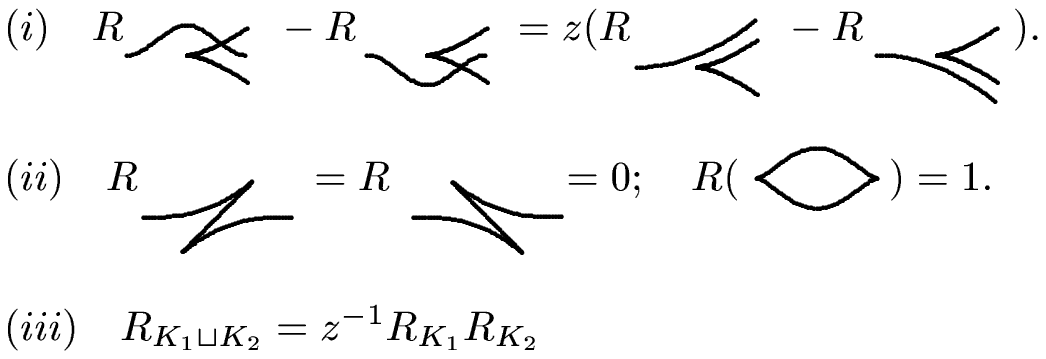}
\end{figure}

\paragraph{Proof for $R_K$.}  Let $L_1, L_2,L_3, L_4$ denote the link diagrams 
appearing from left to right in relation (i).  Divide the rulings 
$\Gamma (L_1)$ into two subsets, those where the visible crossing is switched, 
denoted $A(L_1)$, and those where it is not $B(L_1)$.  Do the same for 
$\Gamma(L_2)$.  There are obvious bijections $B(L_1) \leftrightarrow B(L_2)$ , 
$A(L_1) \leftrightarrow \Gamma(L_3)$, and $A(L_2) \leftrightarrow \Gamma(L_
4)$.  The first preserves the number of switches and the second two decrease 
the number by one.  Hence from (2), \[\begin{array}{rl} R_{L_1} - R_{L_2} &= 
\sum_{\rho \in \Gamma(L_1)}z^{j(\rho)} - \sum_{\rho \in \Gamma(L_2)}z^{j(
\rho)}\\ &= \sum_{\rho \in A(L_1)}z^{j(\rho)} - \sum_{\rho \in A(L_2)}z^{j(
\rho)}\\ &=\sum_{\rho \in \Gamma(L_3)}z^{j(\rho)+1} - \sum_{
\rho \in \Gamma(L_4)}z^{j(\rho)+1}=z(R_{L_3}-R_{L_4}).\end{array}\] The 
relation (ii) is obvious, and the relation (iii) follows from a bijection 
$\theta : \Gamma(K_1) \times \Gamma(K_2) \mapsto \Gamma(K_1\sqcup K_2)$ such 
that $j(\theta(\rho_1,\rho_2)) = j(\rho_1) + j(\rho_2) -1.$  

\paragraph{Proof for $B_K$.}  (i) follows since the corresponding pieces of 
topological diagrams are precisely those in the skein relation for $D_K$ and 
all diagrams have the same number of cusps.  

(ii) follows since we have observed that $B_K$ vanishes when $\beta(K)$ is not 
maximal.  

(iii) follows from the formula $D_{K_1 \sqcup K_2} =\displaystyle{a-a^{-1}
\over z}D_{K_1}D_{K_2})$ \cite{K1}.

To prove that the relations in Lemma 3.1 uniquely characterize a Legendrian 
Isotopy invariant, we realize any Legendrian link as the product of certain 
planar tangles.

\begin{figure}[hbtp]
\centering
\includegraphics[width=5.4in]{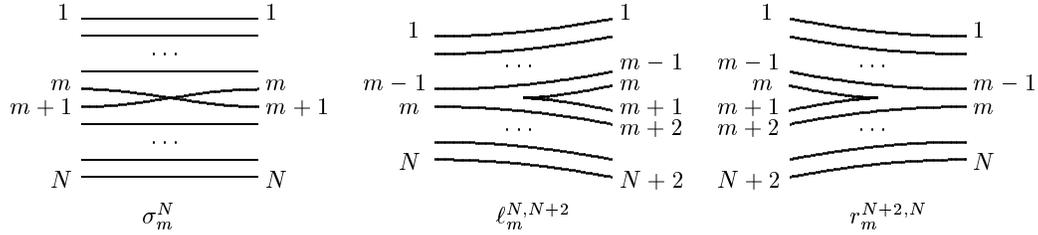}
\caption{Elementary tangles}
\end{figure}

Let $\sigma^N_m$ be the $m$th generator of the $N$ stranded braid group with 
no over or understrand specified at the crossing, $1 \leq m \leq N-1$.  
Let $l^{N,N+2}_m$ be $N$ horizontal lines with a new left cusp appearing 
between the $m-1$th and $m$th strand (increasing the number of strands by 2 in 
the process). Possible values for $m$ are $ 1,\ldots,N,N+1$.  $r^{N+2,N}_m$ 
is its mirror about a vertical line  (see Figure 5.) Compositions are defined 
only when the number of strands agrees, and we will hence omit the upper 
indices from our notation when it will (hopefully) not cause confusion.  
Certainly not every well defined product represents a Legendrian link but 
after a planar isotopy every Legendrian link can be represented by such a 
product.  There are several relations corresponding to Legendrian 
Reidemeister moves and planar isotopy:

Planar isotopy:
\[ \begin{array}{ll} \strut\sigma_{m_1}\sigma_{m_2} = \sigma_{m_2}
\sigma_{m_1},&\mbox{if}\ |m_1-m_2| \geq 2;\\
\strut l_{m_1}^{N, N+2}\sigma_{m_2}^{N+2} = \sigma_{m_2}^N 
l_{m_1}^{N, N+2},&\mbox{if}\ m_1 > m_2+1;\\
\strut l_{m_1}^{N, N+2}\sigma_{m_2}^{N+2} = \sigma_{m_2-2}^N 
l_{m_1}^{N, N+2},&\mbox{if}\ m_2 > m_1+1;\\
\strut\sigma_{m_2}^{N+2}r_{m_1}^{N+2, N} = r_{m_1}^{N+2, N}
\sigma_{m_2}^N, &\mbox{if}\ m_1 > m_2+1;\\
\strut\sigma_{m_2}^{N+2}r_{m_1}^{N+2, N} =  r_{m_1}^{N+2, N}
\sigma_{m_2-2}^N,  &\mbox{if}\ m_2 > m_1+1;\\
\strut l_{m_2}l_{m_1} = l_{m_1-2}l_{m_2}, r_{m_1}r_{m_2} = r_{m_2}r_{m_1-2}, 
&\mbox{if}\ m_1 > m_2+1;\\
\strut r_{m_1}l_{m_2} = l_{m_2} r_{m_1+2}, r_{m_2}l_{m_1} = 
l_{m_1 + 2}r_{m_2},&\mbox{if}\ m_1 \geq m_2.\end{array}\]

Type 1:
$$ l_m\sigma_{m-1}r_m = \id = l_m\sigma_{m+1}r_m$$

Type 2: 
$$ l_{m-1}\sigma_m\sigma_{m-1} = l_m = l_{m+1}\sigma_m\sigma_{m+1}$$

Type 3:
$$ \sigma_{m+1}\sigma_{m}\sigma_{m+1}=\sigma_{m}\sigma_{m+1}\sigma_{m} $$

The skein relations can be realized as \[ \begin{array} {ll} (i)&R_{\dots 
l_{m+1}\sigma_m \ldots} - R_{\ldots l_m \sigma_{m+1}\ldots} = z ( R_{\ldots 
l_{m+1}\ldots} - R_{\ldots l_m\ldots})\strut\\(ii)&R_{\ldots l_m r_{m-1}
\ldots} = R_{\ldots l_m r_{m+1}\ldots} = 0 , R_{\ldots \sigma_ir_i \ldots} = 
R_{\ldots l_i\sigma_i\ldots } =0 , R_{l^{0, 2}_1 r^{2,0}_1}= 1\strut\\ 
(iii)&R_{K_1 \sqcup K_2}=z^{-1}R_{K_1}R_{K_2}\end{array}\]

The second entry of (ii) is implied by the first, but we include it for 
convenience.

\begin{lemma}
By repeated evaluation of the skein relation a formula for 
$R_K$ can be found in terms of the R polynomials of Legendrian links with less 
crossings and links whose values are specified by (ii).
\end{lemma}

\paragraph{Proof.}  Note that since the fronts on the RHS of (i) have less 
crossings then the fronts on the LHS, if the theorem holds for one of the 
fronts on the LHS it must hold for the other.  As a consequence, given a link 
described as a word, $W$, in the above planar tangles, it is enough to show
the following statement. 

(A) By substituting into the above relations corresponding to Legendrian 
isotopy and interchanging $l_{m+1}\sigma_m $ and $l_m \sigma_{m+1}$,  $W$ may 
be reduced to a word with less crossings or to a word whose $R$ polynomial is 
known by (ii). 

We refer to the interchange of $l_{m+1}\sigma_m $ and $l_m \sigma_{m+1}$ as a 
{\sl skein move}.

Statement (A) is proved by nested inductions, the outer being on , $L := $
the number of left cusps of $W$.  The base case is handled entirely by (ii). 
Now, assuming the statement for $L-1$ we prove it for $L$ by the following 
induction.

By looking at the left cusp located farthest to the right in $W$ we can write 
$W = Xl^{N-2,N}_mY$ where $Y$ is a word in the $\sigma_i$ and $r_i$.  Our 
inner induction is on $M :=N + cr(Y)$. The base case $(M = 2)$ is handled by 
(iii) and the outer inductive hypothesis (for then we have a disjoint copy of 
the Legendrian unknot at the end of $W$).

For general $M$, given a word of the form $$Xl_m(\sigma_{m-1}\sigma_{m-2}
\ldots\sigma_{m-N_1})(\sigma_{m+1}\sigma_{m+2}\ldots\sigma_{m+N_2})Y; N_1,N_2
\geq 0$$ we give a procedure depending on the first letter of $Y$ to either 
reduce one of $N$ or $cr((\sigma_{m-1}\sigma_{m-2}\ldots\sigma_{m-N_1})(
\sigma_{m+1}\sigma_{m+2}\ldots\sigma_{m+N_2})Y)$ or increase one of $N_1$ and 
$N_2$.  Since $N_1$ and $N_2$ can only be increased a finite number of times 
this will complete the inductive step.

Note that $(\sigma_{m-1}\sigma_{m-2}\ldots\sigma_{m-N_1})$ commutes with 
$(\sigma_{m+1}\sigma_{m+2}\ldots\sigma_{m+N_2})$ by planar isotopy.

Case 1: $Y = \sigma_iY'$

SubCase 1:  $i< m - N_1 -1$ or $i > m+N_2 +1$.

By planar isotopy $\sigma_i$ commutes with$$l_m(\sigma_{m-1}\sigma_{m-2}
\ldots\sigma_{m-N_1})(\sigma_{m+1}\sigma_{m+2}\ldots\sigma_{m+N_2})$$ so can 
be absorbed into $X$ decreasing $cr(Y)$.

SubCase 2:  $i = m-N_1-1 $ or $i = m+N_2+1$.

Add  $\sigma_i$ at the end of the appropriate parenthesis expression 
increasing $N_1$ or $N_2$.

SubCase 3: $i = m-N_1$ or $i = m+N_2$ but $i \neq m$.

We deal with the first possibility since the second is similar.  Note that 
after a skein move,$$l_m(\sigma_{m-1}\sigma_{m-2}\ldots\sigma_{m-N_1})(
\sigma_{m+1}\sigma_{m+2}\ldots\sigma_{m+N_2})$$becomes $$ l_{m-1}(
\sigma_{m-2}\ldots\sigma_{m-N_1})(\sigma_{m}\sigma_{m+1}\ldots\sigma_{m+
N_2}).$$Applying the skein move $N_1$ times we arrive at\[\begin{array} {l} 
\qquad Xl_{m-N_1}(\sigma_{m-N_1 +1}\ldots\sigma_{m+1}\sigma_{m+2}\ldots
\sigma_{m+N_2})Y\strut\\ =Xl_{m-N_1}(\sigma_{m-N_1 +1}\ldots\sigma_{m+1}
\sigma_{m+2}\ldots\sigma_{m+N_2})\sigma_{m-N_1}Y'\strut\\ =Xl_{m-N_1}(
\sigma_{m-N_1 +1}\sigma_{m-N_1}\sigma_{m-N_1+2}\ldots\sigma_{m+1}\sigma_{m+2}
\ldots\sigma_{m+N_2})Y'\end{array}\]for which a type II Legendrian 
Reidemeister move removes two crossings.

SubCase 4: $ m -N_1 >i > m$  or $m < i < m+N_2$

Again we consider just the first possibility,\[\begin{array} {l} \qquad Xl_m(
\sigma_{m-1}\sigma_{m-2}\ldots\sigma_{m-N_1})(\sigma_{m+1}\sigma_{m+2}\ldots
\sigma_{m+N_2})Y\strut\\ =Xl_m(\sigma_{m-1}\sigma_{m-2}\ldots\sigma_{m-N_1})(
\sigma_{m+1}\sigma_{m+2}\ldots\sigma_{m+N_2})\sigma_iY'\strut\\ =Xl_m(
\sigma_{m-1}\ldots\sigma_i\sigma_{i-1}\sigma_i\ldots\sigma_{m-N_1})(
\sigma_{m+1}\sigma_{m+2}\ldots\sigma_{m+N_2})Y'\strut\\ =Xl_m(\sigma_{m-1}
\ldots\sigma_{i-1}\sigma_{i}\sigma_{i-1}\ldots\sigma_{m-N_1})(\sigma_{m+1}
\sigma_{m+2}\ldots\sigma_{m+N_2})Y'\strut\\ =X\sigma_{i-1}l_m(\sigma_{m-1}
\ldots\sigma_{i}\sigma_{i-1}\ldots\sigma_{m-N_1})(\sigma_{m+1}\sigma_{m+2}
\ldots\sigma_{m+N_2})Y'\end{array}\]decreasing $cr(Y)$ by 1.  The 3rd equality 
is a type 3 Legendrian Reidemeister move.  The rest are planar isotopy.

SubCase 5:  $i = m$

If exactly one of $N_1$ or $N_2$ is $0$ a type 2 Reidemeister can be applied. 
If they are both $0$ the value of the polynomial is $0$ by (ii).  If neither 
are zero we have\[\begin{array}{l}\quad Xl_m(\sigma_{m-1}\sigma_{m-2}\ldots
\sigma_{m-N_1})(\sigma_{m+1}\sigma_{m+2}\ldots\sigma_{m+N_2})Y\strut\\ 
=Xl_m(\sigma_{m-1}\sigma_{m-2}\ldots\sigma_{m-N_1})(\sigma_{m+1}\sigma_{m+2}
\ldots\sigma_{m+N_2})\sigma_mY'\strut\\ =Xl_m\sigma_{m-1}\sigma_{m+1}
\sigma_m(\sigma_{m-2}\ldots\sigma_{m-N_1})(\sigma_{m+2}\sigma_{m+2}\ldots
\sigma_{m+N_2})Y'\strut\\ =Xl_{m-1}\sigma_{m}\sigma_{m+1}\sigma_m(\sigma_{m-2}
\ldots\sigma_{m-N_1})(\sigma_{m+2}\sigma_{m+2}\ldots\sigma_{m+N_2})Y'\strut\\ 
=Xl_{m-1}\sigma_{m+1}\sigma_{m}\sigma_{m+1}(\sigma_{m-2}\ldots\sigma_{m-N_1})
(\sigma_{m+2}\sigma_{m+2}\ldots\sigma_{m+N_2})Y'\strut\\ =X\sigma_{m+1}
l_{m-1}\sigma_{m}\sigma_{m+1}(\sigma_{m-2}\ldots\sigma_{m-N_1})(\sigma_{m+2}
\sigma_{m+2}\ldots\sigma_{m+N_2})Y'\end{array}\]decreasing c(Y) by 1.  
The equalities are assumption, planar isotopy, skein move (i), type 3 
Legendrian Reidemeister move, and planar isotopy respectively.

Case 2. $Y = r_iY'$

Again by vertical symmetry of all relations involved we assume without loss of 
generality that $i \leq m$.

SubCase 1:  $i< m - N_1 -1$

By planar isotopy $r_i$ commutes with$$l_m(\sigma_{m-1}\sigma_{m-2}\ldots
\sigma_{m-N_1})(\sigma_{m+1}\sigma_{m+2}\ldots\sigma_{m+N_2})$$so can be 
absorbed into $X$ decreasing $N$.

SubCase 2:  $i = m-N_1-1 $ 

Applying the skein move $N_1$ times we arrive at $$Xl_{m-N_1}r_{m-N_1-1}
\ldots$$ which is has a zig-zag.

SubCase 3:  $i = m-N_1$

The inclusion of the factor $\sigma_{m-N_1}r_{m-N_1}$ shows that the 
polynomial is 0 by (ii).

SubCase 4: $ m -N_1 >i > m$ 

A type 2 Legendrian Reidemeister move may be applied after a planar isotopy.
\[\begin{array} {l} \quad Xl_m(\sigma_{m-1}\sigma_{m-2}\ldots\sigma_{m-N_1})(
\sigma_{m+1}\sigma_{m+2}\ldots\sigma_{m+N_2})r_iY'\strut\\ =Xl_m(\sigma_{m-1}
\ldots\sigma_i\sigma_{i-1}r_i\ldots\sigma_{m-N_1})(\sigma_{m+1-2}\sigma_{m+2-2}
\ldots\sigma_{m+N_2-2})Y'.\end{array}\]

The presence of $\sigma_i\sigma_{i-1}r_i$ allows two crossings to be removed 
with a type 2 move.

SubCase 5: $i=m$

Again there are 4 cases.  $N_1=N_2=0$ implies we have a disjoint unknot in the 
middle of the link diagram so we apply (iii) and the outer inductive 
hypothesis.  If exactly 1 is zero we can apply a type 1 Reidemeister move to 
remove a crossing.  In the final case we have \[\begin{array}{l} \quad Xl_m(
\sigma_{m-1}\sigma_{m-2}\ldots\sigma_{m-N_1})(\sigma_{m+1}\sigma_{m+2}\ldots
\sigma_{m+N_2})r_mY'\strut\\ =Xl_m\sigma_{m-1}\sigma_{m+1}r_m(\sigma_{m-2}
\ldots\sigma_{m-N_1})(\sigma_{m+2-2}\sigma_{m+2-2}\ldots\sigma_{m+N_2-2})Y'
\strut\\ =Xl_{m-1}\sigma_{m}\sigma_{m+1}r_m(\sigma_{m-2}\ldots\sigma_{m
-N_1})(\sigma_{m+2-2}\sigma_{m+2-2}\ldots\sigma_{m+N_2-2})Y'\end{array}\]where 
we used planar isotopy and the skein move to arrange the sequence 
$\sigma_{m}\sigma_{m+1}r_m$, allowing a type 2 move.

This concludes the proof.

\begin{corollary}
The skein relations in Lemma 3.2 uniquely characterize a Legendrian Isotopy 
invariant.
\end{corollary}

Theorem 3.1 now follows from Lemma 3.2 and Corollary 3.1  

\section{Oriented Rulings and the HOMFLY Polynomial}

\subsection{HOMFLY polynomial}

With an appropriate strengthening of the notion of ruling a very similar 
situation holds with regard to the HOMFLY estimate for $\beta (K)$.  First we 
recall a construction of the HOMFLY polynomial \cite{K2}.  An auxiliary 
polynomial $H(z,a)$ is calculated from oriented diagrams according to the 
skein relations

\begin{figure}[hbtp]
\centering
\includegraphics[width=3.8in]{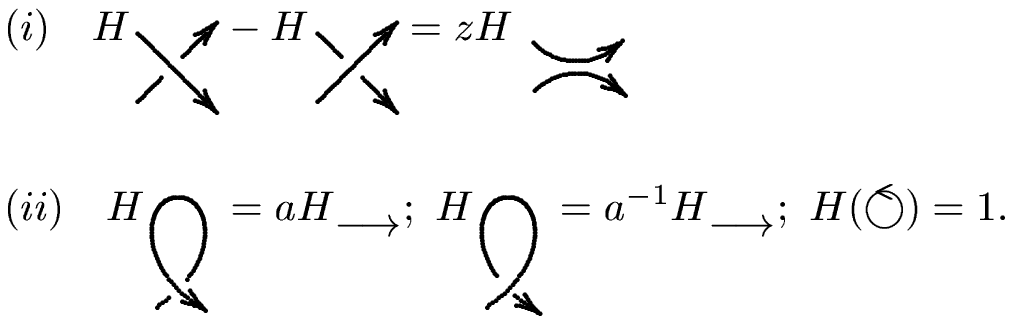}
\end{figure}

\noindent and invariance under regular isotopy.

The HOMFLY polynomial $P_K$ is then defined by the normalization $P_K = 
a^{-w(K)}$.  

\begin{lemma}[\cite{FT}]

$\beta(K) <  -\deg_aP_K$

\end{lemma}

We define an invariant of oriented Legendrian links by $Q_K(z) = C_{c(k) 
-1}(z)$, where $H_K(z,a) = \sum_{n \in {\bf Z}} C_n(z)a^n$.  As in Section 2., 
the non-vanishing of $Q_K(z)$ will be equivalent to the sharpness of Lemma 4.1.

\subsection{Oriented Ruling polynomial}

We call a ruling of a front $K$ {\sl oriented} if all switches are positive 
crossings (in the sense of writhe).  That is, at switches arrows should point 
in the same horizontal direction.  Let $O\Gamma ( K)$ denote the set of 
oriented rulings of a front $K$.  Define the {\sl oriented ruling polynomial} 
of $K$, $$OR_K = \sum_{\rho \in O\Gamma(K)} z^{j(\rho)}$$ $OR_K$ can be seen 
to be a Legendrian isotopy invariant by constructing bijections between 
oriented rulings under Legendrian Reidemeister moves.  

\paragraph{Remark.}  For knots this follows from known results since for knots 
a ruling is oriented if and only if it is 2-graded.  

\begin{theorem}  $OR_K = Q_K$. \end{theorem}

Proof is similar to that of Theorem 1.  Both polynomials are easily seen to 
satisfy the Legendrian skein relations

\begin{figure}[hbtp]
\centering
\includegraphics[width=3.8in]{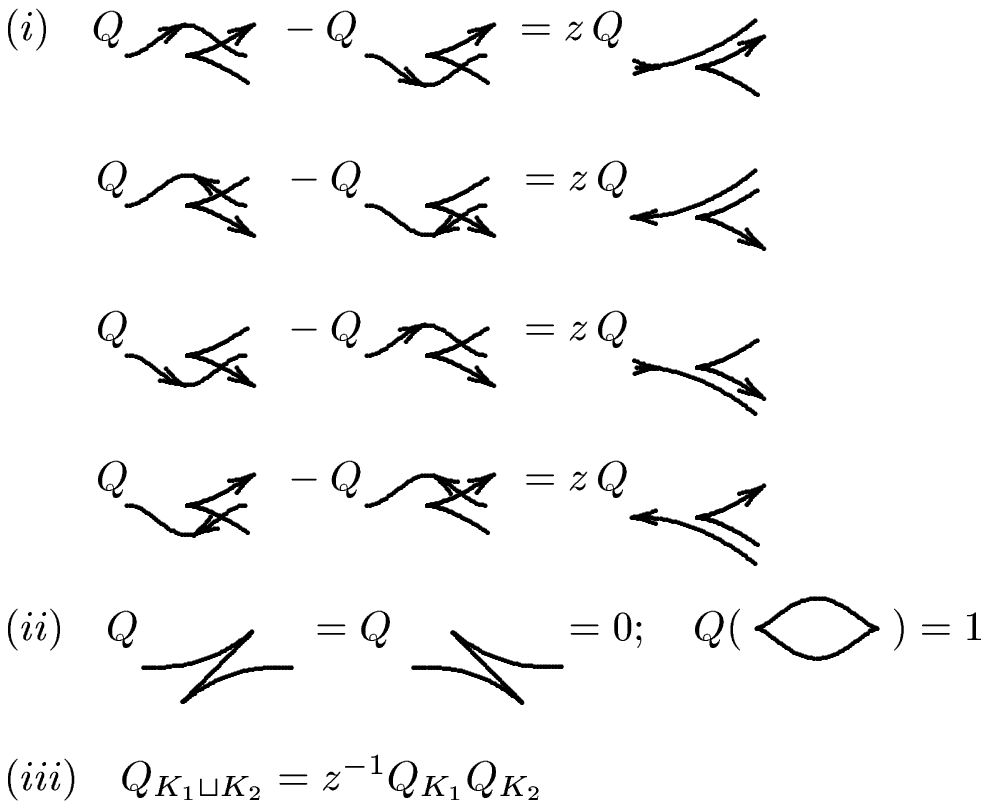}
\end{figure}

\noindent and the same proof shows that these relations uniquely characterize 
a Legendrian isotopy invariant.

\section{Corollaries}

\begin{corollary}  If a Legendrian link $K$ admits a ruling then it maximizes 
$\beta(K)$ within its topological isotopy class.\end{corollary}

\begin{corollary}   The sharpness of the estimate for $\beta (K)$ given by the 
HOMFLY polynomial implies the sharpness of the estimate given by the Kauffman 
polynomial. \end{corollary} 

This is somewhat surprising since it is known that no inequality between the 
two estimates exists in general \cite{Fer}.  As a result we can strengthen the 
estimate given by $P_K$ for some links

\begin{corollary}   If $\deg_aP_K > \deg_aF_K$ then $$\beta(K) <  
-\deg_aP_K-1$$.\end{corollary} 

\begin{corollary}  For an alternating link $L$, $\deg_aP_L \leq \deg_aF_K$.
\end{corollary}

\paragraph{Proof.}  In \cite{Ng1}, Proposition 11 shows that any alternating 
link has a legendrian representative admitting a ruling.  Hence, the estimate 
coming from the Kauffman polynomial is sharp.

This result is conjectured by Ferrand \cite{Fer} using slightly different 
language (\cite{Fer} uses slightly different conventions for the link 
polynomials in question as well as different conventions for front diagrams).

\end{document}